\documentclass{article}
\usepackage{amssymb}
\usepackage{amsfonts}

\input{tcilatex}
\begin{document}

\title{Connected neighborhoods in Cartesian products of solenoids}
\author{Jan P. Boro\'{n}ski, Alejandro Illanes and Emanuel R. M\'{a}rquez}
\maketitle

\begin{abstract}
Given a collection of pairwise co-prime integers $%
m_{1},\ldots ,m_{r}$, greater than $1$, we consider the product $\Sigma =\Sigma _{m_{1}}\times
\cdots \times \Sigma _{m_{r}}$, where each $\Sigma _{m_{i}}$ is the $m_{i}$-adic
solenoid. Answering a question of D. P. Bellamy and J. M. \L ysko, in this
paper we prove that if $M$ is a subcontinuum of $\Sigma $ such that the
projections of $M$ on each $\Sigma _{m_{i}}$ are onto, then for each open
subset $U$ in $\Sigma $ with $M\subset U$, there exists an open connected
subset $V$ of $\Sigma $ such that $M\subset V\subset U$; i.e. any such $M$ is ample in the sense of Prajs and Whittington [10]. This contrasts with the property of Cartesian squares of fixed solenoids  $\Sigma _{m_{i}}\times
 \Sigma _{m_{i}}$, whose diagonals are never ample [1]. 
\end{abstract}

Key words and phrases: solenoids, connected neighborhoods, continuum, fupcon property, product.

Mathematics Subjects Classification: primary 54F15, secondary 54F50.

\bigskip

\begin{center}
\textbf{Introduction}
\end{center}

\bigskip

In the present paper we study Cartesian products of Vietors solenoids. Recall
that a (Vietoris) $n$-adic \textit{solenoid} is a compact and connected space (i.e. \textit{continuum}) given as the inverse limit of circles, with $n$-fold covering maps as bonding maps ($n>1$). Solenoids are \textit{indecomposable}; i.e. they cannot be given as the union of two proper subcontinua. Since solenoids are topological groups they are homogeneous. Solenoids arise in the theory of dynamical systems as suspensions of odometers, and therefore minimal sets of suspension flows. Solenoids were also used by Smale to provide one of the first examples of hyperbolic attractors [11]. Williams [12] generlized these examples to higher dimensional analogues of solenoids, first constructed by McCord [8]. Here we study a condition under which subcontinua of products of solenoids have arbitrarily small connected open neighborhoods, called \textit{full projection implies arbitrary small
	connected open neighborhoods}, or \textit{fupcon} for short. Recall that, given a family of metric continua $\{X_{\alpha
}:\alpha \in J\}$, the product $X=\tprod\nolimits_{\alpha \in J}X_{\alpha }$
has fupcon property, provided that for every subcontinuum $%
M$ and open subset $U$ of $X$ such that $M\subset U$ and $\pi _{\alpha
}(M)=X_{\alpha }$ for each $\alpha \in J$ ($\pi _{\alpha }$ is the $\alpha ^{%
\text{th}}$-projection), there exists an open connected subset $V$ of $X$
such that $M\subset V\subset U$. Clearly, each product of locally connected continua has fupcon property. A subcontinuum $M$ of a continuum $X$ is \textit{ample }provided that for
each open subset $U$ of $X$ with $M\subset U$, there exists a subcontinuum $%
L $ of $X$ such that $M\subset $ int$_{X}(L)\subset L\subset U$. So $X$ is
connected im kleinen at a point $p\in X$ provided that $\{p\}$ is ample. By
[3, Lemma 1], the product $X=\tprod\nolimits_{\alpha \in J}X_{\alpha }$ has
fupcon property provided that each subcontinuum $M$ of $X$ projecting onto
each $X_{\alpha }$ is an ample subset of $X$. Ample subcontinua where introduced in [10] and they have been useful to
improve the understanding of homogeneous continua. It is easy to show that if $M$ is an ample subcontinuum of a continuum $X$,
then the hyperspace $C(X)$ of subcontinua of $X$ (with the Hausdorff metric)
is connected im kleinen at $M$. Thus, if $X$ is a product with fupcon
property, then it is possible to find subcontinua $M$ of $X$ at which $C(X)$
is connected im kleinen. This is something remarkable, since in some of the
examples considered in this area (products of pseudo-arcs, solenoids and
Knaster continua), of products with fupcon property, the properties of local
connectedness are very rare. Below we list some of the known results related to fupcon property.
\begin{enumerate}
	\item any product of Knaster continua has fupcon property [1, Theorem 4.1 and
Observation in p. 230],
\item  any product of pseudo-arcs has fupcon property [1, Theorem 4.4 and
Observation in p. 230],

\item  a product of a solenoid with itself does not have fupcon property [1,
Corollary 4.3],

\item  a product of the pseudo-arc with any product of Knaster continua has
fupcon property [3, Corollary 9],

\item  there exist indecomposable chainable continua whose product does not have
fupcon property [3, Example 11],

\item  if a product of continua has fupcon property, then each factor is a
Kelley continuum [3, Theorem 10],

\item there is a complete characterization of chainable continua $X$ for which
the diagonal in $X\times X$ has arbitrary small connected neighborhoods [4,
Corollary 3.2]

\item  every product of homogeneous continua having the fixed point property has
fupcon property [5, Theorem 2.1],

\item any product of a solenoid and any Knaster continuum has fupcon property
[2, Theorem 1.2] and [5, Theorem 3.1],

\item there exists a Kelley continuum $X$ such that $X\times \lbrack 0,1]$ does
not have fupcon property [5, Example 4.1], and

\item  the product of a chainable Kelley continuum and $[0,1]$ has 
fupcon property [5, Theorem 5.4],

\item  the product of two Smith's nonmetric pseudo-arcs has the fupcon property
[2, Theorem 3.1].

\end{enumerate}
Answering a question of D. P. Bellamy and J. M. \L ysko [1, Question (4)],
in this paper we prove the following theorem.

\bigskip

\textbf{Theorem.} Let $m_{1},\ldots ,m_{r}$ be a finite sequence of
pairwise co-prime integers greater than $1$. For each $i\in \{1,\ldots
,r\}$, let $\Sigma _{m_{i}}$ be the $m_{i}$-adic solenoid. Then $\Sigma =\Sigma
_{m_{1}}\times \cdots \times \Sigma _{m_{r}}$ has the fupcon property.

\bigskip

\begin{center}
\textbf{Covering spaces}
\end{center}

\bigskip

A \textit{mapping} is a continuous function. A \textit{covering space} of a
space $Y$ is a pair $(X,\varphi )$, where $\varphi :X\rightarrow Y$ is a
mapping and for each $y\in Y$, there exists an open neighborhood $U$ of $y$
in $Y$ such that $\varphi ^{-1}(U)$ is a disjoint union of open sets of $X$
each of which is mapped homeomorphically onto $U$ by $\varphi $. Given a covering space $(X,\varphi )$ of the space $Y$ and a mapping $\psi
:Z\rightarrow Y$, a \textit{lifting }of $\psi $ to $X$ is a mapping $%
h:Z\rightarrow X$ such that $\psi =\varphi \circ h$. In the case that there
exists such an $h$, we say that $\psi $ \textit{can be lifted}. By $S^{1}$ we denote the unit circle (centered at the origin) in the
Euclidean plane. By $\mathbb{R}$ we denote the real line. We consider the
covering space $(\mathbb{R},e)$ of $S^{1}$, where $e$ is the exponential
mapping $e:\mathbb{R}\rightarrow S^{1}$ defined by $e(t)=(\cos
(2\pi t),\sin (2\pi t))$. We consider also covering spaces of the form $%
((S^{1})^{r},g)$ and $((S^{1})^{r},g^{n})$, where $(S^{1})^{r}$ is the
Cartesian product $S^{1}\times \cdots \times S^{1}$ (with $r$ factors), $%
g:(S^{1})^{r}\rightarrow (S^{1})^{r}$ is a mapping of the form $%
g(z_{1},\ldots ,z_{r})=(z_{1}^{k_{1}},\ldots ,z_{r}^{k_{r}})$, where $%
k_{1},\ldots ,k_{r}\in \mathbb{N}$ and $g^{n}$ is the composition $g\circ
\cdots \circ g$ ($n$ times). The following lemma is an immediate consequence of Theorem 5.1 of [9].

\bigskip

\textbf{Lemma 1. }Let $U$ be a connected locally arcwise connected space, $%
u_{0}\in U$ and $g:U\rightarrow S^{1}$ a mapping that cannot be lifted to a mapping $\tilde{g}: U\to \mathbb{R}$, with $g=e\circ\tilde{g}$. Then there exists a loop $\gamma
:[0,1]\rightarrow (U,u_{0})$ such that $g\circ \gamma :[0,1]\rightarrow
(S^{1},g(u_{0}))$ is a non-trivial loop.

\bigskip

\begin{center}
\textbf{Loops in }$(S^{1})^{r}$
\end{center}

\bigskip

From now on we fix a finite sequence of integers $m_{1},\ldots ,m_{r}$ such
that if $i\neq j$, then $m_{i}$ and $m_{j}$ are relatively prime, and $%
m_{i}\geq 2$ for each $i$. We also fix the mapping $f:(S^{1})^{r}\rightarrow
(S^{1})^{r}$ given by

\begin{center}
$f(z_{1},\ldots ,z_{r})=(z_{1}^{m_{1}},\ldots ,z_{r}^{m_{r}})$.
\end{center}
We denote by $1^{-}$ the point in $(S^{1})^{r}$ given by $%
1^{-}=((1,0),\ldots ,(1,0))$. For each $i\in \{1,\ldots ,r\}$, we denote by $\rho
_{i}:(S^{1})^{r}\rightarrow S^{1}$ the $i^{\text{th}}$-projection.

\bigskip

\begin{center}
\textbf{A general construction}
\end{center}
\bigskip

Given $k_{1},\ldots ,k_{r}\in \mathbb{N}$, consider the mapping $%
g:(S^{1})^{r}\rightarrow (S^{1})^{r}$ given by $g(z_{1},\ldots
,z_{r})=(z_{1}^{k_{1}},\ldots ,z_{r}^{k_{r}})$. Then $((S^{1})^{r},g)$ is a
covering space of the space $(S^{1})^{r}$. Given a loop $\gamma :[0,1]\rightarrow ((S^{1})^{r},1^{-})$, define $\gamma
^{\ast }:[0,\infty )\rightarrow (S^{1})^{r}$ by
\begin{center}
$\gamma ^{\ast }(t)=\gamma (t-i+1)$, if $t\in \lbrack i-1,i]$.
\end{center}
Since the fundamental group of $[0,\infty )$ is trivial, by Theorem 5.1 of
[9], there exists a unique lifting $\gamma (g):[0,\infty )\rightarrow
(S^{1})^{r}$ of the mapping $\gamma ^{\ast }$ (in the covering space $%
((S^{1})^{r},g)$) such that $\gamma (g)(0)=1^{-}$. Then $\gamma (g)$ satisfies the following properties.

\smallskip

(a) For each $t\in \lbrack 0,\infty )$, $g(\gamma (g)(t))=\gamma (t-k+1)$,
if $t\in \lbrack k-1,k]$,
\smallskip

(b) $g(\func{Im}(\gamma (g)))\subset \func{Im}(\gamma )$.

\bigskip

\noindent
For each $k\in \mathbb{N}$, $g(\gamma (g)(k))=\gamma (1)=1^{-}$. Then $%
\gamma (g)(\mathbb{N})\subset g^{-1}(1^{-})$. Thus, $\gamma (g)(\mathbb{%
N\cup }\{0\})\,$\ is finite. Moreover, $\gamma (g)([k-1,k])=$ the image of
the lifting of $\gamma $ that starts at the point $\gamma (g)(k-1)$. Thus:
\smallskip

(c) $\func{Im}\gamma (g)$ is a subcontinuum of $(S^{1})^{r}$.

\bigskip

\noindent Notice that for each $n\in \mathbb{N}$, we can apply the above construction
to the mapping $g=f^{n}$. Then for each loop $\gamma :[0,1]\rightarrow
((S^{1})^{r},1^{-})$ and for each $n\in \mathbb{N}$, we can define the
mapping $\gamma ^{(n)}:[0,\infty )\rightarrow (S^{1})^{r}$ by

\begin{center}
$\gamma ^{(n)}=\gamma (f^{n})$.
\end{center}
Given $s_{0}=(s_{1},\ldots ,s_{r})\in (\mathbb{Z}\smallsetminus \{0\})^{r}$
and $n\in \mathbb{N}\cup \{0\}$, let $\sigma (s_{0},n):[0,\infty
)\rightarrow (S^{1})^{r}$ be given by

\begin{center}
$\sigma (s_{0},n)(t)=(e(\frac{s_{1}t}{m_{1}^{n}}),\ldots ,e(\frac{s_{r}t}{%
m_{r}^{n}}))$.
\end{center}
For each $s\in \mathbb{Z}\smallsetminus \{0\}$, let $\lambda
_{s}:[0,1]\rightarrow S^{1}$ be given by

\begin{center}
$\lambda _{s}(t)=e(st)$.
\end{center}

\bigskip

\textbf{Lemma 2. }Let $n,k\in \mathbb{N}$, $s_{0}=(s_{1},\ldots ,s_{r})\in (%
\mathbb{Z}\smallsetminus \{0\})^{r}$ and $\gamma ,\zeta :[0,1]\rightarrow
((S^{1})^{r},1^{-})$ loops. Then

(d) $f\circ \gamma ^{(n+1)}=\gamma ^{(n)}$,

(e) $(\sigma (s_{0},0)|_{[0,1]})^{(n)}=(\sigma
(s_{0},0)|_{[0,1]})(f^{n})=\sigma (s_{0},n)$,

(f) If $\gamma $ and $\zeta $ are homotopic, then $\gamma ^{(n)}|_{[0,k]}$
and $\zeta ^{(n)}|_{[0,k]}$ are homotopic, and $\gamma ^{(n)}(k)=\zeta
^{(n)}(k)$,

(g) If $\gamma :[0,1]\rightarrow ((S^{1})^{r},1^{-})$ is a loop such that
for each $i\in \{1,\ldots ,r\}$, $\rho _{i}\circ \gamma $ is homotopic to $%
\lambda _{s_{i}}$, then $\gamma $ is homotopic to $\sigma (s_{0},0)|_{[0,1]}$%
.

\bigskip

\textbf{Proof. }(d) Given $t\in \lbrack 0,\infty )$, suppose that $t\in
\lbrack i-1,i]$ with $i\in \mathbb{N}$, then $f^{n}(f\circ \gamma
^{(n+1)})(t)=f^{n+1}(\gamma ^{(n+1)}(t))=f^{n+1}(\gamma (f^{n+1})(t))=\gamma
(t-i+1)=\gamma ^{\ast }(t)$ and $(f\circ \gamma ^{(n+1)})(0)=f(\gamma
(f^{n+1})(0))=f(1^{-})=1^{-}$. Thus, $f\circ \gamma ^{(n+1)}$ is the unique
lifting of $\gamma ^{\ast }$ (in the covering space $((S^{1})^{r},f^{n})$)
such that $(f\circ \gamma ^{(n+1)})(0)=1^{-}$. Hence, $f\circ \gamma
^{(n+1)}=\gamma ^{(n)}$. 

(e) Given $t\in \lbrack 0,\infty )$, suppose that $t\in \lbrack i-1,i]$ with 
$i\in \mathbb{N}$. By the properties of the exponential mapping, we have that$\ f^{n}(\sigma (s_{0},n)(t))=f^{n}(e(\frac{%
s_{1}t}{m_{1}^{n}}),\ldots ,e(\frac{s_{r}t}{m_{r}^{n}}))=(e(s_{1}t),\ldots
,e(s_{r}t))=(e(s_{1}(t-(i-1))),\ldots ,e(s_{r}(t-(i-1)))=((\sigma
(s_{0},0)_{[0,1]})^{\ast })(t)$.

Then the mapping $\sigma (s_{0},n)$ is a lifting of the mapping $(\sigma
(s_{0},0)|_{[0,1]})^{\ast }$ (in the covering space $((S^{1})^{r},f^{n})$)
such that $\sigma (s_{0},n)(0)=1^{-}$. So $(\sigma
(s_{0},0)|_{[0,1]})(f^{n})=\sigma (s_{0},n)$.

(f) Notice that $\gamma ^{\ast }|_{[0,k]}$ and $\zeta ^{\ast }|_{[0,k]}$ are
homotopic loops in $((S^{1})^{r},1^{-})$. Notice also that $\gamma
^{(n)}|_{[0,k]}$ and $\zeta ^{(n)}|_{[0,k]}$ are liftings of $\gamma ^{\ast
}|_{[0,k]}$ and $\zeta ^{\ast }|_{[0,k]}$ (in the covering space $%
((S^{1})^{r},f^{n})$), respectively such that $\gamma ^{(n)}(0)=1^{-}=\zeta
^{(n)}(0)$. By Lemma 3.3 in [9], $\gamma ^{(n)}|_{[0,k]}$ and $\zeta
^{(n)}|_{[0,k]}$ are homotopic and $\gamma ^{(n)}(k)=\zeta ^{(n)}(k)$.

(g) Is immediate.

\bigskip

\textbf{Theorem 3. }Let $s_{0}=(s_{1},\ldots ,s_{r})\in (\mathbb{Z}%
\smallsetminus \{0\})^{r}$. Then there exists $N\in \mathbb{N}$ such that
for each $n\geq N$ and every loop $\gamma :[0,1]\rightarrow
((S^{1})^{r},1^{-})$ satisfying that for each $i\in \{1,\ldots ,r\}$, $\rho
_{i}\circ \gamma $ is homotopic to $\lambda _{s_{i}}$, we have that $f^{-1}(%
\func{Im}(\gamma ^{(n)}))=\func{Im}(\gamma ^{(n+1)})$ and $f^{-1}(\func{Im}%
(\gamma ^{(n)}))$ is a subcontinuum of $(S^{1})^{r}$.

\bigskip

\textbf{Proof. }Let $\alpha _{1},\ldots ,\alpha _{r}\in \mathbb{N}\cup \{0\}$
be such that for each $i\in \{1,\ldots ,r\}$, $s_{i}=m_{i}^{\alpha
_{i}}q_{i} $, where $q_{i}$ and $m_{i}$ are relative prime. Let $\alpha
=\max \{\alpha _{1},\ldots ,\alpha _{r}\}$, $N=m_{1}^{\alpha }\cdots
m_{r}^{\alpha }$ and $n\geq N$. Let $\gamma :[0,1]\rightarrow (S^{1})^{r}$ be a loop such that for each $%
i\in \{1,\ldots ,r\}$, $\rho _{i}\circ \gamma $ is homotopic to $\lambda
_{s_{i}}$.

\bigskip
\noindent
By Lemma 2 (d), $f\circ \gamma ^{(n+1)}=\gamma ^{(n)}$. Thus, $\func{Im}%
\gamma ^{(n+1)}\subset f^{-1}(\func{Im}(\gamma ^{(n)}))$.

\bigskip
\noindent
By Lemma 2 (g), $\gamma $ is homotopic to $\sigma (s_{0},0)|_{[0,1]}$. By
Lemma 2, (f) and (f), for each $k\in \mathbb{N}$, $\gamma
^{(n+1)}(k)=(\sigma (s_{0},0)|_{[0,1]})^{(n+1)}(k)=\sigma (s_{0},n+1)(k)$.
Since $\gamma ^{(n+1)}(0)=1^{-}=\sigma (s_{0},n+1)(0)$, we conclude that $%
\gamma ^{(n+1)}(k)=\sigma (s_{0},n+1)(k)$ for each $k\in \mathbb{N}\cup
\{0\} $.

\bigskip
\noindent
Notice that $f^{-1}(1^{-})=\{(e(\frac{j_{1}}{m_{1}}),\ldots ,e(\frac{j_{r}}{%
m_{r}})):$ for each $i\in \{1,\ldots ,r\}$, $0\leq j_{i}<m_{i}\}$. We claim that $f^{-1}(1^{-})\subset \gamma ^{(n+1)}(\mathbb{N})$. For
proving this, take $i\in \{1,\ldots ,r\}$ and $j_{i}\in \{0,\ldots
,m_{i}-1\} $. Set $j^{-}=(j_{1},\ldots ,j_{r})$.

\bigskip
\noindent
For each $i\in \{1,\ldots ,r\}$, let $\beta _{i}=n-\alpha _{i}$. Set $%
u=m_{1}^{\beta _{1}}\cdots m_{r}^{\beta _{r}}$ and $u_{i}=\frac{u}{%
m_{i}^{\beta _{i}}}$. Then $\frac{us_{i}}{m_{i}^{n+1}}=\frac{um_{i}^{\alpha
_{i}}q_{i}}{m_{i}^{n+1}}=\frac{u_{i}q_{i}}{m_{i}}$. Notice that $u_{i}q_{i}$
and $m_{i}$ are relative prime. By the Chinese Remainder Theorem, there exists $x\in \mathbb{N}$ such that
for each $i\in \{1,\ldots ,r\}$,

\begin{center}
$u_{i}q_{i}x\equiv j_{i}(\func{mod}m_{i})$.
\end{center}
Then $e(\frac{us_{i}x}{m_{i}^{n+1}})=e(\frac{u_{i}q_{i}x}{m_{i}})=e(\frac{%
j_{i}}{m_{i}})$ for each $i\in \{1,\ldots ,r\}$. Set $k=ux$. Then $k\in \mathbb{N}$ and

\begin{center}
$\gamma ^{(n+1)}(k)=\sigma (s_{0},n+1)(k)=\sigma (s_{0},n+1)(ux)=(e(\frac{%
s_{1}ux}{m_{1}^{n+1}}),\ldots ,e(\frac{s_{r}ux}{m_{r}^{n+1}}))=(e(\frac{j_{1}%
}{m_{1}}),\ldots ,e(\frac{j_{r}}{m_{r}}))$.
\end{center}
We have shown that $f^{-1}(1^{-})\subset \gamma ^{(n+1)}(\mathbb{N})$.

\bigskip
\noindent
Take a point $v\in f^{-1}(\func{Im}(\gamma ^{(n)}))$. Then $f(v)=\gamma
^{(n)}(t_{0})$ for some $t_{0}\in \lbrack 0,\infty )$. Consider the path $\gamma ^{(n)}|_{[0,t_{0}]}$ in $(S^{1})^{r}$. Since $%
((S^{1})^{r},f)$ is a covering space of $(S^{1})^{r}$, there exists a
lifting $\zeta :[0,t_{0}]\rightarrow (S^{1})^{r}$ of $\gamma
^{(n)}|_{[0,t_{0}]}$ such that $\zeta (t_{0})=v$. Since $f\circ \zeta
=\gamma ^{(n)}|_{[0,t_{0}]}$, we have $f(\zeta (0))=\gamma ^{(n)}(0)=1^{-}$,
so $\zeta (0)\in f^{-1}(1^{-})\subset \gamma ^{(n+1)}(\mathbb{N})$. Then
there exists $k\in \mathbb{N}$ such that $\zeta (0)=\gamma ^{(n+1)}(k)$.

\smallskip
\noindent
Consider the mapping $\omega :[0,k+t_{0}]\rightarrow (S^{1})^{r}$ given by:

\begin{center}
$\omega (t)=\left\{ 
\begin{array}{cc}
\gamma ^{(n+1)}(t)\text{,} & \text{if }t\in \lbrack 0,k]\text{,} \\ 
\zeta (t-k)\text{,} & \text{if }t\in \lbrack k,k+t_{0}]\text{.}%
\end{array}%
\right. $
\end{center}
Given $t\in \lbrack 0,k+t_{0}]$, if $t\leq k$, by the definition of $\gamma
^{(n+1)}$, $f^{n+1}(\omega (t))=f^{n+1}(\gamma ^{(n+1)}(t))=f^{n+1}((\gamma
(f^{n+1}))(t))=\gamma (t-j+1)$, if $t\in \lbrack j-1,j]$ ($j\in \mathbb{N}$%
). If $k\leq t$, by the definition of $\gamma ^{(n)}$, $f^{n+1}(\omega
(t))=f^{n+1}(\zeta (t-k))=f^{n}(\gamma ^{(n)}(t-k))=\gamma (t-k-j+1)$, if $%
t-k\in \lbrack j-1,j]$. That is, $f^{n+1}(\omega (t))=\gamma (t-k-j+1)$, if $%
t\in \lbrack k+j-1,k+j]$.

\smallskip
\noindent
We have shown that for each $t\in \lbrack 0,k+t_{0}]$, $f^{n+1}(\omega
(t))=\gamma (t-i+1)$, if $t\in \lbrack i-1,i]$ ($i\in \mathbb{N}$). Since $%
\gamma ^{(n+1)}$ also satisfies that for each $t\in \lbrack 0,k+t_{0}]$, $%
f^{n+1}(\gamma ^{(n+1)}(t))=\gamma (t-i+1)$, if $t\in \lbrack i-1,i]$ and $%
\omega (0)=\gamma ^{(n+1)}(0)$, we conclude that $\omega =\gamma
^{(n+1)}|_{[0,k+t_{0}]}$.

\smallskip
\noindent
In particular, $\gamma ^{(n+1)}(k+t_{0})=\omega (k+t_{0})=\zeta (t_{0})=v$.
Therefore, $v\in \func{Im}(\gamma ^{(n+1)})$. This completes the proof that $%
\func{Im}(\gamma ^{(n+1)})=f^{-1}(\func{Im}(\gamma ^{(n)}))$.

\smallskip
\noindent
Finally, by (c), $\func{Im}(\gamma ^{(n+1)})$ is a subcontinuum of $%
(S^{1})^{r}$. $\blacksquare $

\bigskip

\begin{center}
\textbf{Fupcon property}\ 
\end{center}

\bigskip

For each $i\in \{1,\ldots ,r\}$, we consider the solenoid $\Sigma _{m_{i}}$
which is the inverse limit of the inverse sequence $\{S^{1},f_{i}\}$, where $%
f_{i}:S^{1}\rightarrow S^{1}$ is defined as before. That is $%
f_{i}(z)=z^{m_{i}}$. We consider the Cartesian product

\begin{center}
$\Sigma =\Sigma _{m_{1}}\times \cdots \times \Sigma _{m_{r}}$
\end{center}
For each $i\in \{1,\ldots ,r\}$, let $\eta _{i}:\Sigma \rightarrow \Sigma
_{m_{i}}$ be the $i^{\text{th}}$-projection. We identify $\Sigma $ with $\lim_{\leftarrow }((S^{1})^{r},f)$, where $f$ is
defined as before. That is

\begin{center}
$f(z)=(f_{1}(z),\ldots ,f_{r}(z))$.
\end{center}
For each $n\in \mathbb{N}$, let $\pi _{n}:\Sigma \rightarrow (S^{1})^{r}$ be
the projection given by:

\begin{center}
$\pi _{n}(z_{1},z_{2},\ldots )=z_{n}$, where each $z_{i}\in (S^{1})^{r}$.
\end{center}

\bigskip

\textbf{Lemma 4. }Let $M$ be a subcontinuum of $\Sigma $ such that the point 
$p_{0}=(1^{-},1^{-},\ldots )$ belongs to $M$ and $\eta _{i}(M)=\Sigma
_{m_{i}}$ for each $i\in \{1,\ldots ,r\}$. Then for each $n\in \mathbb{N}$
and each $i\in \{1,\ldots ,r\}$, $(\rho _{i}\circ \pi
_{n})|_{M}:M\rightarrow S^{1}$ cannot be lifted (in the covering space $(%
\mathbb{R},e)$).

\bigskip

\textbf{Proof. }Suppose to the contrary that $(\rho _{i}\circ \pi _{n})|_{M}$
can be lifted for some $n\in \mathbb{N}$ and $i\in \{1,\ldots ,r\}$. Then
there exists a mapping $h:M\rightarrow \mathbb{R}$ such that $e\circ h=(\rho
_{i}\circ \pi _{n})|_{M}$. Since $(\rho _{i}\circ \pi _{n})(p_{0})=(1,0)$,
we assume that $h(p_{0})=0$.

\smallskip
\noindent
Suppose that $i=1$, the other cases are similar.
For each $j\in \{1,\ldots ,m_{1}\}$, define $g_{j}:M\rightarrow \mathbb{R}$
by $g_{j}(p)=\frac{h(p)+j-1}{m_{1}}$.
Let $K_{j}=\{p\in M:\rho _{1}(\pi _{n+1}(p))=e(g_{j}(p))\}$. Then $K_{j}$ is
closed in $M$.

\smallskip
\noindent
Given $p=(p_{1},p_{2},\ldots )\in M$, where for each $k\in \mathbb{N}$, $%
p_{k}=(p_{1}^{(k)},\ldots ,p_{r}^{(k)})$, we have $f(p_{n+1})=p_{n}$, so $%
(\rho _{1}(\pi _{n+1}(p)))^{m_{1}}=(p_{1}^{(n+1)})^{m_{1}}=p_{1}^{(n)}=(\rho
_{1}\circ \pi _{n})(p)=e(h(p))=(e(\frac{h(p)}{m_{1}}))^{m_{1}}$. This
implies that $\frac{\rho _{1}(\pi _{n+1}(p))}{e(\frac{h(p)}{m_{1}})}$ is an $%
m_{1}$-root of the unity in the complex plane $\mathbb{C}$. Thus, there
exists $j\in \{1,\ldots ,m_{1}\}$ such that $\rho _{1}(\pi _{n+1}(p))=e(%
\frac{h(p)}{m_{1}})e(\frac{j-1}{m_{1}})=e(\frac{h(p)+j-1}{m_{1}})$. Hence, $%
p\in K_{j}$. We have shown that $M=K_{1}\cup \cdots \cup K_{m_{1}}$.

\smallskip
\noindent
Take $j,l\in \{1,\ldots ,m_{1}\}$ such that $j\neq l$. If there exists $p\in
K_{j}\cap K_{l}$, then $e(g_{j}(p))=\rho _{1}(\pi _{n+1}(p))=e(g_{l}(p))$.
This implies that $\frac{h(p)+j-1}{m_{1}}=\frac{h(p)+l-1}{m_{1}}+k$ for some 
$k\in \mathbb{Z}$. Then $h(p)+j-1=h(p)+l-1+km_{1}$, and $j-l=km_{1}$. This
contradicts the fact that $j\neq l$ and $j,l\in \{1,\ldots ,m_{1}\}$. We
have shown that $K_{1},\ldots ,K_{m_{1}}$ are pairwise disjoint.

\smallskip
\noindent
Since $\rho _{1}(\pi _{n+1}(p_{0}))=(1,0)=e(\frac{h(p_{0})}{m_{1}})$, we
obtain that $K_{1}\neq \emptyset $. The connectedness of $M$ implies that $K_{1}=M$. Thus, for each $p\in M\,$, $%
\rho _{1}(\pi _{n+1}(p))=e(\frac{h(p)}{m_{1}})$. Let $h_{1}=\frac{1}{m_{1}}h$. Then $h_{1}$ is a lifting of $\rho _{1}\circ
\pi _{n+1}$ such that $h_{1}(p_0)=\frac{h(p_{0})}{m_{1}}=0$. If we repeat the argument
with $h_{1}$ instead of $h$, we can obtain that for each $p\in M$, $\rho
_{1}(\pi _{n+2}(p))=e(\frac{h(p)}{m_{1}^{2}})$.

\smallskip
\noindent
Proceeding this way we can obtain that for each $a\in \mathbb{N}$, $\rho
_{1}(\pi _{n+a}(p))=e(\frac{h(p)}{m_{1}^{a}})$.

\smallskip
\noindent
Since $M$ is compact, there exists $a_{0}\in \mathbb{N}$ such that diameter$%
(\{\frac{h(p)}{m_{1}^{a_{0}}}:p\in M\})<1$. Thus, $\rho _{1}(\pi
_{n+a_{0}}(M))=\{e(\frac{h(p)}{m_{1}^{a_{0}}}):p\in M\}\neq S^{1}$. Therefore, there
exists a point $w_{n+a_{0}}\in S^{1}\smallsetminus \rho _{1}(\pi _{n+a_{0}}(M))$. Take $w_{0}\in \Sigma _{m_{1}}$ be any point of the form

\begin{center}
$w_{0}=(w_{1},\ldots ,w_{n+a-1},w_{n+a},w_{n+a+1},\ldots )$.
\end{center}

Then no point $p$ in $M$ satisfies $\eta _{1}(p)=w_{0}$. This contradiction
finishes the proof of the lemma. $\blacksquare $

\bigskip

\textbf{Lemma 5. }Let $M$ be a subcontinuum of $\Sigma $ such that the point 
$p_{0}=(1^{-},1^{-},\ldots )$ belongs to $M$ and $\eta _{i}(M)=\Sigma
_{m_{i}}$ for each $i\in \{1,\ldots ,r\}$. Then for each $n\in \mathbb{N}$,
and each open connected subset $U$ of $(S^{1})^{r}$ with $\pi _{n}(M)\subset
U$, there exists a loop $\gamma $ in $(U,1^{-})$ such that for each $i\in
\{1,\ldots ,r\}$, the loop $\rho _{i}\circ \gamma :[0,1]\rightarrow S^{1}$
cannot be lifted in the covering space $(\mathbb{R},e)$.

\bigskip

\textbf{Proof. }By Lemma 4, $(\rho _{1}\circ \pi _{n})|_{M}:M\rightarrow
S^{1}$ cannot be lifted. This implies that $\rho _{1}|_{\pi _{n}(M)}:\pi
_{n}(M)\rightarrow S^{1}$ cannot be lifted, and then $\rho
_{1}|_{U}:U\rightarrow S^{1}$ cannot be lifted.

\smallskip
\noindent
By Lemma 1, there exists a loop $\gamma _{1}:[0,1]\rightarrow (U,1^{-})$
such that $\rho _{1}\circ \gamma _{1}:[0,1]\rightarrow (S^{1},(1,0))$ is a
non-trivial loop. Thus, there exists $s\in \mathbb{Z}\smallsetminus \{0\}$
such that $\rho _{1}\circ \gamma _{1}$ is homotopic to $\lambda _{s}$.

\smallskip
\noindent
Suppose, inductively, that $1\leq k<r$ and we have constructed a loop $%
\gamma _{k}:[0,1]\rightarrow (U,1^{-})$ such that for each $i\in \{1,\ldots
,k\}$, the mapping $\rho _{i}\circ \gamma _{k}:[0,1]\rightarrow (S^{1},(1,0))
$ is homotopic to $\lambda _{s_{i}}$ for some $s_{i}\in \mathbb{Z}%
\smallsetminus \{0\}$. Suppose that $\rho _{k+1}\circ \gamma _{k}$ (which is
a loop in $(S^{1},(1,0))$) is homotopic to $\lambda _{s_{k+1}}$ for some $%
s_{k+1}\in \mathbb{Z}$.

\smallskip
\noindent
Using Lemma 1 again, we have that there exists a loop $\zeta
:[0,1]\rightarrow (U,1^{-})$ such that $\rho _{k+1}\circ \zeta
:[0,1]\rightarrow (S^{1},(1,0))$ is homotopic to $\lambda _{s_{0}}$ for some 
$s_{0}\in \mathbb{Z}\smallsetminus \{0\}$.

\smallskip
\noindent
For each $i\in \{1,\ldots ,k\}$, let $t_{i}\in \mathbb{Z}$ be such that $%
\rho _{i}\circ \zeta $ is homotopic to $\lambda _{t_{i}}$.
Let $l\in \mathbb{N}$ be such that $\max \{\left\vert s_{1}\right\vert
,\ldots ,\left\vert s_{k}\right\vert ,\left\vert s_{k+1}\right\vert \}<l\leq
l\left\vert s_{0}\right\vert $.

\smallskip
\noindent
Let $\gamma _{k+1}=\zeta \ast \cdots \ast \zeta \ast \gamma
_{k}:[0,1]\rightarrow U$, where $\zeta $ appears $l$ times in $%
\zeta \ast \cdots \ast \zeta $. Then $\gamma _{k+1}$ is a loop in $%
U$.

\smallskip
\noindent
Given $i\in \{1,\ldots ,k+1\}$, $\rho _{i}\circ \gamma _{k+1}=(\rho
_{i}\circ \zeta )\ast \cdots \ast (\rho _{i}\circ \zeta )\ast (\rho
_{i}\circ \gamma _{k})$ is homotopic to $\lambda _{u_{i}}$ for some $%
u_{i}\in \mathbb{Z}$.

\smallskip
\noindent
If $i\leq k$, then $u_{i}=lt_{i}+s_{i}$. In the case that $t_{i}=0$, then $%
u_{i}=s_{i}\neq 0$, and in the case that $t_{i}\neq 0$, then $l\left\vert
t_{i}\right\vert \geq l>\left\vert s_{i}\right\vert $, so $u_{i}\neq 0$.

\smallskip
\noindent
If $i=k+1$, then $u_{i}=ls_{0}+s_{k+1}\neq 0$.We have shown that $u_{i}\neq 0$ for every $i\in \{1,\ldots ,k+1\}$. This
finishes the induction and the proof of the lemma. $\blacksquare $

\bigskip

\textbf{Theorem 6. }$\Sigma =\Sigma _{m_{1}}\times \cdots \times \Sigma
_{m_{r}}$ has the fupcon property.

\bigskip

\textbf{Proof. }Let $M$ be a subcontinuum of $\Sigma $ such that for each $%
i\in \{1,\ldots ,r\}$, $\eta _{i}(M)=\Sigma _{m_{i}}$. Let $\varepsilon >0$. We are going to show that there exists a connected open subset $V$ of $%
\Sigma $ such that $M\subset V$ and $V$ is contained in the $\varepsilon $%
-neighborhood of $M$ in the space $\Sigma $.

\smallskip
\noindent
Since $\Sigma $ is homogeneous, we may assume that $p_{0}=(1^{-},1^{-},%
\ldots )\in M$. Let $N_{0}\in \mathbb{N}$ be such that $\frac{1}{2^{N_{0}}}<\frac{%
\varepsilon }{2}$.

\smallskip
\noindent
By the uniform continuity of $f$, there exists $\delta >0$ such that $\delta
<\varepsilon $ and the following implication holds: if $i\in \{1,\ldots
,N_{0}\}$ and $z,w\in (S^{1})^{r}$ are such that $\left\vert z-w\right\vert
<\delta $, then $\left\vert f^{i}(z)-f^{i}(w)\right\vert <\frac{\varepsilon 
}{2}$.

\smallskip
\noindent
Let $W$ be the $\delta $-neighborhood of $\pi _{N_{0}}(M)$ in the space $%
(S^{1})^{r}$.

\smallskip
\noindent
By Lemma 5, there exists a loop $\gamma _{0}$ in $(W,1^{-})$ such that for
each $i\in \{1,\ldots ,r\}$, the loop $\rho _{i}\circ \gamma
_{0}:[0,1]\rightarrow (S^{1},(1,0))$ cannot be lifted in $(\mathbb{R},e)$.
Therefore, there exists $s_{i}\in \mathbb{Z}\smallsetminus \{0\}$ such that $%
\rho _{i}\circ \gamma _{0}$ is homotopic to $\lambda _{s_{i}}$.

\smallskip
\noindent
Let $N_{1}\in \mathbb{N}$ be as in Theorem 3 such that $N_{1}>N_{0}$.
Let $Y$ be the component of $f^{-N_{1}}(W)$ that contains $\pi
_{N_{0}+N_{1}}(M)$. Then $Y$ is a connected open subset of $(S^{1})^{r}$ and 
$1^{-}\in Y$.

\smallskip
\noindent
Let $V=\pi _{N_{0}+N_{1}}^{-1}(Y)$. Then $V$ is open in $\Sigma $ and $%
M\subset V$.
We are going to show that $V$ is connected.

\smallskip
\noindent
Take $v=(v_{1},v_{2},\ldots )\in V$. Then $v_{N_{0}+N_{1}}\in Y$. Let $\eta :[0,1]\rightarrow (S^{1})^{r}$ be a lifting (in the covering space 
$((S^{1})^{r},f^{N_{1}})$) of the loop $\gamma _{0}$ such that $\eta
(0)=1^{-}$. Then $f^{N_{1}}(\eta (1))=\gamma _{0}(1)=1^{-}$. Since $%
f^{N_{1}}(\func{Im}(\eta ))\subset \func{Im}(\gamma _{0})\subset W$, we have 
$\func{Im}(\eta )\subset f^{-N_{1}}(W)$, and since $1^{-}\in \func{Im}(\eta
)\cap \pi _{N_{0}+N_{1}}(M)$, we conclude that $\func{Im}(\eta )\subset Y$.

\smallskip
\noindent
Since $Y$ is arcwise connected, there exists a path $\mu :[0,1]\rightarrow Y$
such that $\mu (0)=\eta (1)$ and $\mu (1)=v_{N_{0}+N_{1}}$.

\smallskip
\noindent
Let $\omega $ be the mapping $\eta \ast \mu \ast \mu ^{-1}:[0,1]\rightarrow
(S^{1})^{r}$ and let $\gamma :[0,1]\rightarrow (S^{1})^{r}$ be the
loop $f^{N_{1}}\circ \omega :[0,1]\rightarrow ((S^{1})^{r},1^{-})$.
Notice that $\func{Im}(\gamma )\subset W$, $\gamma $ is homotopic to $%
f^{N_{1}}\circ \eta =\gamma _{0}$. By the choice of $N_{1}$, for each $n\geq
N_{1}$, $f^{-1}(\func{Im}(\gamma ^{(n)}))$ is a subcontinuum of $(S^{1})^{r}$%
.

\smallskip
\noindent
Notice that $\omega $ is a lifting of the loop $\gamma $ (in the covering
space $((S^{1})^{r},f^{N_{1}})$) satisfying $\omega (0)=\eta (0)=1^{-}$. On
the other hand, $\gamma ^{(N_{1})}|_{[0,1]}$ is a lifting (in $%
((S^{1})^{r},f^{N_{1}})$) of $\gamma $ with $\gamma
(f^{N_{1}})(0)=1^{-}$. Thus, $\omega =\gamma ^{(N_{1})}|_{[0,1]}$.
Therefore, $\func{Im}(\omega )\subset \func{Im}(\gamma ^{(N_{1})})$. In
particular, $v_{N_{0}+N_{1}}\in \func{Im}(\gamma ^{(N_{1})})$.

\smallskip
\noindent
By (b) $f^{N_{1}}(\func{Im}(\gamma (f^{N_{1}})))\subset \func{Im}(\gamma )$, so $\func{Im}(\gamma ^{(N_{1})})=\func{Im}(\gamma
(f^{N_{1}}))\subset f^{-N_{1}}(W)$. Moreover, since $(\gamma
^{(N_{1})})(0)=1^{-}$, we have that $\func{Im}(\gamma ^{(N_{1})})\subset Y$.

\smallskip
\noindent
Consider the sequence $\{L_{n}\}_{n=1}^{\infty }$ of subcontinua of $%
(S^{1})^{r}$ defined in the following way.
\begin{center}
$L_{1}=f^{N_{0}-1}(\func{Im}(\gamma ))$, $L_{2}=f^{N_{0}-2}(\func{Im}(\gamma
))$, \ldots\ , $L_{N_{0}-1}=f(\func{Im}(\gamma ))$, $L_{N_{0}}=\func{Im}%
(\gamma )$,

$L_{N_{0}+1}=\func{Im}(\gamma ^{(1)})$, $L_{N_{0}+2}=\func{Im}(\gamma
^{(2)}) $, $L_{N_{0}+3}=\func{Im}(\gamma ^{(3)})$, \ldots\ , $%
L_{N_{0}+N_{1}}=\func{Im}(\gamma ^{(N_{1})})$,

$L_{N_{0}+N_{1}+1}=f^{-1}(\func{Im}(\gamma ^{(N_{1})}))$, $%
L_{N_{0}+N_{1}+2}=f^{-1}(\func{Im}(\gamma ^{(N_{1}+1)}))$, $%
L_{N_{0}+N_{1}+3}=f^{-1}(\func{Im}(\gamma ^{(N_{1}+2)}))$, \ldots
\end{center}

\bigskip

\textbf{Claim.}

(1) $\{(L_{n},f)\}_{n=1}^{\infty }$ is an inverse sequence of continua,

(2) the inverse limit $L_{0}$ of the sequence $\{(L_{n},f)\}_{n=1}^{\infty }$
is a subcontinuum of $\Sigma $ containing $v$, contained in $V$ and
intersecting $M$.

(3) $V$ is contained in the $\varepsilon $-neighborhood of $M$ in the space $%
\Sigma $.

\bigskip

\textbf{Proof of (1).} By (c), $L_{1},\ldots ,L_{N_{0}+N_{1}}$ are
subcontinua of $(S^{1})^{r}$. The choice of $N_{1}$ implies that $%
L_{N_{0}+N_{1}+1}$, $L_{N_{0}+N_{1}+2}$, \ldots\ are also continua.

\smallskip
\noindent
Clearly, for each $1<n\leq N_{0}$, $f(L_{n})=L_{n-1}$. By (b), $%
f(L_{N_{0}+1})\subset L_{N_{0}}$.

\smallskip
\noindent
By Lemma 2 (e), for each $n\in \mathbb{N}$, $$f(\func{Im}(\gamma
^{(n+1)}))\subset \func{Im}(\gamma ^{(n)})$$ 
and 
$$f(f^{-1}(\func{Im}(\gamma
^{(n+1)})))\subset \func{Im}(\gamma ^{(n+1)})\subset f^{-1}(\func{Im}(\gamma
^{(n)})).$$ So, for each $n\in \{N_{0}+1,\ldots ,N_{0}+N_{1}-1\}\cup
\{N_{0}+N_{1}+1,N_{0}+N_{1}+2,\ldots \}$, $f(L_{n+1})\subset L_{n}$.
Clearly, $f(L_{N_{0}+N_{1}+1})\subset L_{N_{0}+N_{1}}$. We have shown that
for each $n\in \mathbb{N}$, $f(L_{n+1})\subset L_{n}$. This completes the
proof of (1).

\bigskip

\textbf{(2). }By the choice of $N_{1}$, $L_{N_{0}+N_{1}+2}=f^{-1}(\func{Im}(\gamma
^{(N_{1}+1)}))=f^{-2}(\func{Im}(\gamma ^{(N_{1})}))$, $%
L_{N_{0}+N_{1}+3}=f^{-3}(\func{Im}(\gamma ^{(N_{1})}))$, $%
L_{N_{0}+N_{1}+4}=f^{-4}(\func{Im}(\gamma ^{(N_{1})}))$,\ldots 

\smallskip
\noindent
Since $v_{N_{0}+N_{1}}\in \func{Im}(\gamma ^{(N_{1})})$ and $v\in \Sigma $, $%
v_{N_{0}+N_{1}+1}\in f^{-1}(\func{Im}(\gamma ^{(N_{1})}))=L_{N_{0}+N_{1}+1}$%
, $v_{N_{0}+N_{1}+2}\in f^{-2}(\func{Im}(\gamma
^{(N_{1})}))=L_{N_{0}+N_{1}+2}$ and so on. Moreover, by (1), $$%
v_{N_{0}+N_{1}-1}\in L_{N_{0}+N_{1}-1},$$ $$v_{N_{0}+N_{1}-2}\in
L_{N_{0}+N_{1}-2}$$ and so on. Therefore, $v\in L_{0}$.

\smallskip
\noindent
By (b), $f^{N_{1}}(\func{Im}(\gamma ^{(N_{1})}))=f^{N_{1}}(\func{Im}(\gamma
(f^{N_{1}})))\subset \func{Im}(\gamma )\subset W$. So, $\func{Im}(\gamma
^{(N_{1})})\subset f^{-N_{1}}(W)$. Recall that $1^{-}\in \func{Im}(\gamma
^{(N_{1})})$. This implies that $L_{N_{0}+N_{1}}=\func{Im}(\gamma
^{(N_{1})})\subset Y$. Thus, $L_{0}\subset \pi _{N_{0}+N_{1}}^{-1}(Y)=V$.

\smallskip
\noindent
Since $f(1^{-})=1^{-}$ and $\gamma (0)=f^{N_{1}}(\omega (0))=1^{-}$, it
follows that $1^{-}\in L_{n}$ for all $n\in \mathbb{N}$. Therefore, $%
(1^{-},1^{-},\ldots )\in L_{0}\cap M$.

\bigskip

\textbf{(3). }Take an element $w=(w_{1},w_{2},\ldots )\in V$. Then $%
w_{N_{0}+N_{1}}\in Y$ and $w_{N_{0}}=f^{N_{1}}(w_{N_{0}+N_{1}})\in W$.

\smallskip
\noindent
Since $W$ is the $\delta $-neighborhood of $\pi _{N_{0}}(M)$ in $(S^{1})^{r}$%
, there is $z=(z_{1},z_{2},\ldots )\in M$ such that $\left\vert
w_{N_{0}}-z_{N_{0}}\right\vert <\delta <\varepsilon $. This implies that for
each $i\in \{1,\ldots ,N_{0}-1\}$, $\left\vert
f^{i}(w_{N_{0}})-f^{i}(z_{N_{0}})\right\vert <\frac{\varepsilon }{2}$. Thus,
for each $i\in \{1,\ldots ,N_{0}\}$, $\left\vert w_{i}-z_{i}\right\vert <%
\frac{\varepsilon }{2}$. Since $\frac{1}{2^{N_{0}}}<\frac{\varepsilon }{2}$,
we conclude that the distance in $\Sigma $ from $w$ to $z$ is less than $%
\varepsilon $.

\bigskip
\noindent
Property (2) implies that $V$ is connected. By property (3) we conclude
that $\Sigma $ has the fupcon property. $\blacksquare $

\begin{center}
\textbf{Open problems}
\end{center}
\bigskip

\noindent
Below we list some related open problems. 
\begin{itemize} 
	
	\item[(Q1)] Does the product of two chainable Kelley continua have the fupcon
property [3, Question 12]?

\item[(Q2)] Suppose $X$ and $Y$ are $1$-dimensional continua such that $X\times Y$
has the fupcon property, does the product $X_{P}\times Y_{P}$ have the
fupcon property [2, Question 4]?, $X_{p}$ is a $1$-dimensional continuum
that admits a continuous decomposition into pseudo-arcs, and whose
decomposition space is homeomorphic to $X$ [7]?

\item[(Q3)] Does the product of $[0,1]$ and the pseudo-circle have the fupcon
property [2, Question 5]?

\item[(Q4)] Does the product of a pseudo-arc and pseudo-circle have the fupcon
property [2, Question 6]?

\item[(Q5)] Does the product of two pseudo-circles have the fupcon property [2,
Question 7]?

\item[(Q6)] Cartesian products of which matchbox manifolds have the fupcon property?
\end{itemize}

\bigskip

\begin{center}
\textbf{Acknowledgments}
\end{center}

\bigskip
\noindent
The research in this paper was carried out during the 12th Research Workshop
in Hyperspaces and Continuum Theory held in the city of Quer\'{e}taro, M\'{e}%
xico, during June, 2018. The authors would like to thank Jorge M. Mart\'{\i}%
nez-Montejano who joined the discussion on the topic of this paper.

\bigskip
\noindent
This paper was partially supported by the project "Teor\'ia de Continuos, Hiperespacios y Sistemas Din\'amicos III" (IN106319) of PAPIIT, DGAPA, UNAM. The first author was also supported by University of Ostrava grant lRP201824 "Complex topological structures" and the NPU II project LQ1602 IT4Innovations excellence in science.

\bigskip

\begin{center}
\textbf{References}
\end{center}

\bigskip
\noindent
[1] D. P. Bellamy and J. M. \L ysko, \textit{Connected open neighborhoods of
subcontinua of product continua with indecomposable factors}, Topology Proc.
44 (2014), 223-231.

\noindent
[2] J. P. Boro\'{n}ski, D. R. Prier, M. Smith, \textit{Ample continua in
Cartesian products of continua}, Topology Appl. 238 (2018), 54-58.

\noindent
[3] A. Illanes, \textit{Connected open neighborhoods in products}, Acta
Math. Hungar. 148 (1) (2016), 73-82.

\noindent
[4] A. Illanes, \textit{Small connected neighborhoods containing the
diagonal of a product}, Topology Appl. 230 (2017), 506-516.

\noindent
[5] A. Illanes, J. M. Mart\'{\i}nez-Montejano and K. Villarreal, \textit{%
Connected neighborhoods in products}, Topology Appl. 241 (2018), 172-184.

\noindent
[6] A. Illanes and S. B. Nadler, Jr., \textit{Hyperspaces, Fundamentals and
Recent Advances}, Monographs and Textbooks in Pure and Applied Math., Vol.
216, Marcel Dekker, Inc., New York, Basel, 1999.

\noindent
[7] W. Lewis, \textit{Continuous curves of pseudo-arcs}, Houston J. Math. 11
(1985), 91-99.

\noindent
[8] C. McCord, \textit{Inverse limit sequences with covering maps}, Trans. Amer. Math. Soc. 114 (1965),197–209

\noindent
[9] W. S. Massey, \textit{Algebraic Topology, An Introduction}, Graduate
Texts in Mathematics, v. 56, Springer-Verlag, New York, Heidelberg, Berlin,
1967.

\noindent
[10] J.R. Prajs and K. Whittington, \textit{Filament sets, aposyndesis, and
the decomposition theorem of Jones}, Trans. Amer. Math. Soc. 359 (2007),
5991-6000.

\noindent
[11] S. Smale, {\em Differentiable dynamical systems,} Bull. Amer. Math. Soc., 73 (1967), 747--817.

\noindent
[12] R.F. Williams, {\em Expanding attractors}, Publ. Math. I.H.E.S. 43 (1974), 169–203
\bigskip

Addresses

\bigskip

Boro\'{n}ski: National Supercomputing Centre IT4Innovations, Division of the
University of Ostrava, Instituto for Research and Applications of Fussy
Modeling, 30. Dubna 22, 701 03 Ostrava. Czech Republic, and Faculty of
Applied Mathematics, AGH University of Science and Technology, al.
Mickiewicza 30, 30-059 Krak\'{o}w, Poland.

\bigskip

Illanes: Instituto de Matem\'{a}ticas, Universidad Nacional Aut\'{o}noma de M%
\'{e}xico, Circuito Exterior, Cd. Universitaria, M\'{e}xico, D.F., 04510, M%
\'{e}xico.

\bigskip 

M\'{a}rquez:  Departamento de Matem\'aticas, Facultad de Ciencias, Universidad Nacional Aut\'onoma de M\'exico, Circuito Exterior, Cd. Universitaria, M\'exico, D.F., 04510, M\'exico

\bigskip 

Boro\'{n}ski: jan.boronski@osu.cz

Illanes: illanes@matem.unam.mx

M\'{a}rquez: emanuelrmarquez@outlook.com

\end{document}